\documentclass[a4paper, 11pt]{article}
\usepackage[utf8]{inputenc}
\usepackage{amsmath,amssymb,exscale}
\usepackage{a4wide}
\usepackage[english]{babel}

\begin{document}

\newcommand{\Z}{\ensuremath{\mathbb{Z}}}
\newcommand{\C}{\ensuremath{\mathbb{C}}}
\newcommand{\N}{\ensuremath{\mathbb{N}}}
\newcommand{\Q}{\ensuremath{\mathbb{Q}}}
\newcommand{\R}{\ensuremath{\mathbb{R}}}
\newcommand{\G}{\ensuremath{\mathbb{G}}}
\renewcommand{\P}{\ensuremath{\mathbb{P}}}
\renewcommand{\H}{\ensuremath{\mathrm{H}}}
\newcommand{\M}{\ensuremath{\mathcal{M}}}
\newcommand{\mC}{\ensuremath{\mathcal{C}}}
\newcommand{\mD}{\ensuremath{\mathcal{D}}}
\newcommand{\mE}{\ensuremath{\mathcal{E}}}
\newcommand{\mF}{\ensuremath{\mathcal{F}}}
\newcommand{\mG}{\ensuremath{\mathcal{G}}}
\newcommand{\mH}{\ensuremath{\mathcal{H}}}
\newcommand{\mM}{\ensuremath{\mathcal{M}}}
\newcommand{\mT}{\ensuremath{\mathcal{T}}}
\newcommand{\mU}{\ensuremath{\mathcal{U}}}
\newcommand{\U}{\ensuremath{\mathcal{U}}}
\renewcommand{\_}{\underline}
\renewcommand{\Im}{\ensuremath{\mathrm{Im}}}
\newcommand{\tr}{\ensuremath{\mathrm{tr}}}
\newcommand{\id}{\ensuremath{\mathrm{id}}}
\newcommand{\Spec}{\ensuremath{\mathrm{Spec\,}}}
\newcommand{\rk}{\ensuremath{\mathrm{rk\,}}}
\renewcommand{\c}{\ensuremath{\mathrm{c}}}
\newcommand{\Pic}{\ensuremath{\mathrm{Pic}}}
\newcommand{\mfo}{\ensuremath{\mathfrak{o}}}

\newcommand{\ch}{\ensuremath{\mathrm{ch}}}
\newcommand{\td}{\ensuremath{\mathrm{td}}}
\newcommand{\Coh}{\ensuremath{\mathrm{Coh}}}
\newcommand{\ext}{\ensuremath{\mathrm{ext}}}
\newcommand{\Ext}{\ensuremath{\mathrm{Ext}}}
\renewcommand{\hom}{\ensuremath{\mathrm{hom}}}
\newcommand{\emor}{\ensuremath{\mathrm{end}}}
\newcommand{\Hom}{\ensuremath{\mathrm{Hom}}}
\newcommand{\Ob}{\ensuremath{\mathrm{Ob}}}
\newcommand{\Mor}{\ensuremath{\mathrm{Mor}}}
\newcommand{\End}{\ensuremath{\mathrm{End}}}
\newcommand{\Num}{\ensuremath{\mathrm{Num}}}
\renewcommand{\O}{\ensuremath{\mathcal{O}}}
\newcommand{\Gl}{\ensuremath{\mathrm{Gl}}}
\newcommand{\Sl}{\ensuremath{\mathrm{Sl}}}
\newcommand{\PGl}{\ensuremath{\mathrm{PGl}}}
\newcommand{\Aut}{\ensuremath{\mathrm{Aut}}}
\newcommand{\PAut}{\ensuremath{\mathrm{PAut}}}
\newcommand{\ds}{\ensuremath{/\!\!/}}
\newcommand{\Supp}{\ensuremath{\mathrm{Supp}}}
\newcommand{\Jac}{\ensuremath{\mathrm{Jac}}}
\newcommand{\Hilb}{\ensuremath{\mathrm{Hilb}}}
\newcommand{\Quot}{\ensuremath{\mathrm{Quot}}}
\newcommand{\Drap}{\ensuremath{\mathrm{Drap}}}
\newcommand{\Kern}{\ensuremath{\mathrm{kern}}}
\newcommand{\codim}{\ensuremath{\mathrm{codim}}}
\newcommand{\NS}{\ensuremath{\mathrm{NS}}}
\newcommand{\NSXQ}{\ensuremath{\mathrm{NS}(X)_\Q}}
\newcommand{\Amp}{\ensuremath{\mathrm{Amp}}}
\newcommand{\Hev}{\ensuremath{\Lambda(X)}}
\newcommand{\Hevp}{\ensuremath{\Lambda_+(X)}}
\newcommand{\sm}{\ensuremath{\setminus\{0\}}}
\newcommand{\lel}{\ensuremath{\;(\le)\;}}
\newcommand{\lle}{\ensuremath{\;(<)\;}}
\newcommand{\geg}{\ensuremath{\;(\ge)\;}}
\newcommand{\lelo}{\ensuremath{\;(\le_0)\;}}
\newcommand{\gego}{\ensuremath{\;(\ge_0)\;}}
\newcommand{\mand}{\ensuremath{\quad\mathrm{and}\quad}}
\newcommand{\mwith}{\ensuremath{\quad\mathrm{with}\quad}}
\newcommand{\lcm}{\ensuremath{\mathrm{lcm}}}
\newcommand{\sgn}{\ensuremath{\mathrm{sgn}}}

\newcommand{\Proof}{\noindent\textit{{Proof. }}}
\newcommand{\ProofOf}{\noindent\textit{{Proof of }}}
\newcommand{\qed}{\hfill $\Box$}

\newtheorem{Prop}{Proposition}[section]
\newtheorem{Def}[Prop]{Definition}
\newtheorem{Cor}[Prop]{Corollary}
\newtheorem{Claim}[Prop]{Claim}
\newtheorem{Th}[Prop]{Theorem}
\newtheorem{Lemma}[Prop]{Lemma}
\newtheorem{Obs}[Prop]{Observation}
\newtheorem{Qu}[Prop]{Question}
\newtheorem{Exa}[Prop]{Example}
\newtheorem{Rem}[Prop]{Remark}

\renewcommand{\labelenumi}{\arabic{enumi}.\ }
\renewcommand{\labelenumii}{\alph{enumii})\ }

\selectlanguage{english}
\hyphenation{pa-ra-me-tri-sing co-re-pre-sents se-mi-sta-ble ar-gu-ments ma-ni-folds}

\title{On moduli spaces of sheaves on K3 or abelian surfaces}
\author{Markus Zowislok\footnote{Institute of mathematics, Johannes Gutenberg-Universität Mainz, 55099 Mainz}}
\date{}
\maketitle

\begin{abstract}
We investigate the moduli spaces of one- and two-dimensional sheaves on projective K3 and abelian surfaces that are semistable with respect to a nongeneral ample divisor with regard to the symplectic resolvability.
We can exclude the existence of new examples of projective irreducible symplectic manifolds lying birationally over components of the moduli spaces of one-dimensional semistable sheaves on K3 surfaces, and over components of many of the moduli spaces of two-dimensional sheaves on K3 surfaces, in particular, of those for rank two sheaves.
\end{abstract}

\section{Introduction}
How are moduli spaces of one- and two-dimensional sheaves on a projective K3 or abelian surface $X$ that are semistable with respect to an ample divisor $H$ on $X$ related when $H$ varies in the ample cone? 
Is there a symplectic resolution of the moduli space if $H$ is nongeneral?

For the second question, we follow the idea of constructing a projective $\Q$-factorial symplectic terminalisation $\tilde M\to M$ of a component $M$ of the moduli space, i.e.\ a symplectic $\Q$-factorial projective variety $\tilde M$ with at most terminal singularities together with a projective birational morphism $f: \tilde M\to M$. The existence of such a morphism yields the following facts:
\begin{enumerate}
\item If $\tilde M$ can be chosen to be an irreducible symplectic manifold then $\tilde M$ is unique up to deformation by a result of Huybrechts \cite{Huy}.
\item If $\tilde M$ is singular then $M$ admits no projective symplectic resolution by \cite{Na1} corollary 1.
\end{enumerate}
We denote the moduli space of sheaves on $X$ with Mukai vector $v\in\Hev:=\N_0\oplus\NS(X)\oplus\Z\subset H^{2*}(X,\Z)$ that are semistable with respect to an ample divisor $H$ on $X$ by $M_H(v)$ and the open subscheme of stable sheaves by $M^s_H(v)$. $M^s_H(v)$ is nonsingular, each connected component has dimension $2+v^2$ and it carries a symplectic form due to Mukai \cite{Mu}.
We start with investigating the possible components and show that we can reduce to considering components containing stable sheaves.

We treat the one-dimensional case first.
Let $v=(0,v_1,v_2)\in\Hev$ with $v_1\neq 0$, $v_2\neq 0$, $H$ an ample divisor in a $v$-chamber $K$ and $H'$ another ample divisor in the closure $\overline K$ of $K$ in the ample cone.
Then one has
\begin{center}
$H'$-stable $\Rightarrow$ $H$-stable $\Rightarrow$ $H$-semistable $\Rightarrow$ $H'$-semistable
\end{center}
for one-dimensional sheaves $F$ with $v(F)=v$.
Our main result is the extension of what is known for a general ample divisor:\\

\begin{samepage}
{\noindent \bf Theorem.\ (\ref{MHng0}) \it
Let $v=(0,v_1,v_2)\in\Hev$ with $v_1\neq 0$, $v_2\neq 0$ and $v^2\ge 0$, and $H$ an ample divisor on $X$.
Assume that $M^s_H(v)$ is nonempty.
\begin{enumerate}
\item Let $v$ be primitive or $v^2=8$.
Then there is a projective symplectic resolution
$M\to \overline{M^s_H(v)}$.
If $H$ is not $v$-general then $M$ can be chosen to be a symplectic resolution of $M_A(v)$, where $A$ is a $v$-general ample divisor in a chamber touching $H$.
\item Let $v$ be not primitive and $v^2\neq 8$.
Then there is a singular locally factorial (and therefore $\Q$-factorial) projective symplectic terminalisation 
of $\overline{M^s_H(v)}\,.$
\end{enumerate}
}
\end{samepage}

\noindent 
Together with our discussion above and the decomposition of components of the moduli space containing no stable sheaves this implies that
if for $v=(0,v_1,v_2)\in\Hev$ with $v_1\neq 0$, $H$ an ample divisor on $X$ and $M$ an irreducible component of $M_{H}(v)$
there is a projective symplectic resolution $\tilde M\to M$ with $\tilde M$ an irreducible symplectic manifold
then it is deformation equivalent to a symplectic resolution of some $M_A(w)$, where $w\in\Hev$ and $A$ is some $w$-general ample divisor.\\

For the rest of the article we focus on the two-dimensional case, which is more complicated, as semistable sheaves in general become unstable if the divisor is moved onto a wall.
In \cite{MW} the authors construct a moduli space for twisted semistable sheaves of fixed Chern character on a surface, which can be used to establish a connection between the moduli spaces for varying ample divisor.
They show under certain conditions an equivalence of twisted semistability and semistability with an extra condition involving a second ample divisor $A$.
We will call the latter one $(H,A)$-semistability, which we generalise to the context of a projective scheme in the appendix.
Here $H$ and $A$ are two ample line bundles on the given projective scheme.
The definition immediately yields the observation 
\begin{center}
$H$-stable $\Rightarrow$ $(H,A)$-stable $\Rightarrow$ $(H,A)$-semistable $\Rightarrow$ $H$-semistable,
\end{center}
and thus the needed morphisms between the corresponding moduli spaces.
The construction of the moduli space $M_{H,A}(P)$ of $(H,A)$-semistable sheaves with Hilbert polynomial $P$ with respect to $H$ is also given in the appendix. We generalise the one given in \cite{HL} by using two different ample line bundles $H$ and $A$, the first one in order to make the considered sheaves globally generated, and the second one in order to get the linearised line bundle, similarly to the surface case in \cite{MW}. The reason why we redo the construction is because we need more properties than developed in \cite{MW}. Moreover, although we only use the surface case, we want to remark that this restriction is unnecessary for the construction of the moduli space.

The moduli space $M_{H,A}(v)$ of $(H,A)$-semistable sheaves with Mukai vector $v$ and a $v$-general ample divisor $A$ is a good candidate for a suitable terminalisation as it shares many properties with $M_A(v)$ - 
assuming the existence of certain stable sheaves, which can be ensured by a numerical condition on $v$, see the main text.

Altogether we can extend the results of \cite{KLS} for a general ample divisor $H$ as follows:\\

{\noindent \bf Theorem.\ \it
Let $v=(v_0,v_1,v_2)\in\Hevp:=\N\oplus\NS(X)\oplus\Z$ primitive with $v^2\ge 0$, $m\in\N$ and $H$ an ample divisor on $X$, and assume that $M^s_H(mv)$ is nonempty. 
\begin{enumerate}
\item Let $m=1$ or $(mv)^2=8$.
Then there is a projective symplectic resolution
$M\to \overline{M^s_H(mv)}$.
If $H$ is not $mv$-general then $M$ can be chosen to be a symplectic resolution of $M_{H,A}(mv)$, where $A$ is an $mv$-general ample divisor.
\item Let $m\ge 2$ and $(mv)^2\neq 8$.
If $H$ is $mv$-general or $v_0=1$ or $v^2>\varphi(v_0)$ with $\varphi$ as in theorem \ref{MAHBstab}
then there is a singular $\Q$-factorial projective symplectic terminalisation 
of $\overline{M^s_H(mv)}\,,$
and in particular, there is no projective symplectic resolution of $\overline{M^s_H(mv)}$.\\
\end{enumerate}
}

\Proof
The results for an $mv$-general ample divisor $H$ are well-known, see \cite{KLS}. Let now $H$ be a not $mv$-general ample divisor.
By proposition \ref{MHngPTerm} it is enough to find a $\Q$-factorial projective symplectic terminalisation
$M\to M_{H,A}(v)\,,$ where $A$ is any $mv$-general ample divisor and $M$ is nonsingular in case 1 and singular in case 2. Choose an ample divisor $A$ in an $mv$-chamber touching $H$.
Theorem \ref{MAHgen} yields the claim of item 1 and reduces the claim of item 2 to the question of the existence of an $(H,A)$-stable sheaf with Mukai vector $v$.
If $v_0=1$ then $M_A(v)\neq\emptyset$ (see e.g.\ \cite{KLS}), and any sheaf in this space is also $(H,A)$-stable.
If $v^2>\varphi(v_0)$ then there is an $(H,A)$-stable sheaf with Mukai vector $v$ by theorem \ref{MAHBstab}.
\qed\\

\noindent We have only partial results on the deformation classes of the symplectic manifolds given in part 1 of the theorem:\\

{\noindent \bf Corollary \ref{DefClass}.\ \it
Let $X$ be a K3 surface, $v=(v_0,v_1,v_2)\in\Hevp$ a primitive Mukai vector, $H$ an ample divisor and $A$ a $v$-general ample divisor.
Assume that $v^2> \varphi(v_0)$.
Then $M_{H,A}(v)$ is deformation equivalent to $\Hilb^{\frac{v^2}2+1}(X)$.\\
}

\noindent Remark that our results imply that no new examples of a projective irreducible symplectic manifold arise from moduli spaces of semistable sheaves on projective K3 surfaces of rank two.

\vspace{1cm}
{
\small
\noindent \textit{Acknowledgements.}
The author would like to express his gratitude to his doctoral advisor Manfred Lehn for the optimal support, his suggestions and his encouragement.
He thanks all colleagues for stimulating discussions, in particular Arvid Perego, Christian Lehn, Tanja Becker, Heinrich Hartmann and Oliver Goert\-sches.
Moreover, he also thanks Christoph Sorger for an unpublished article and 
Samuel Boissi\`ere and Olivier Serman for a preprint and corresponding announced results.
Finally, the author thanks the referee for various suggestions for improvements.
Most part of this work was supported by the Sonderforschungsbereich/Transregio 45 ``Periods, moduli spaces and arithmetic of algebraic varieties'' of the Deutsche Forschungsgemeinschaft (German Research Foundation).
This article is partly contained in my PhD thesis \cite{ich}.
}

\section{Notation and conventions}
\subsection{Symplectic varieties}
Following \cite{Beau} a symplectic variety $X$ is a normal variety together with a (holomorphically) symplectic form $\omega$ on the nonsingular locus $U$ of $X$ such that there is a resolution of singularities $f: \tilde X\to X$ for which the pullback $\left(f|_{f^{-1}(U)}\right)^*\omega$ extends to a holomorphic 2-form on $\tilde X$.
One can show that if $X$ is a symplectic variety and $f: \tilde X\to X$ is any resolution of singularities then the pullback of $\omega$ by the induced isomorphism extends to a holomorphic 2-form on $\tilde X$, see \cite{Beau} section 1.
An irreducible symplectic manifold is a simply connected compact Kähler manifold with a holomorphically symplectic form that generates $H^0(X,\Omega_X^2)$. 
\begin{Def}
Let $X$ be a scheme. A nonsingular symplectic variety $\tilde X$ together with a proper birational morphism $f: \tilde X\to X$ is called a symplectic resolution.
\end{Def}
Note that we do not require $f$ to be an isomorphism over the nonsingular locus,
but for a projective symplectic resolution of a projective normal variety this condition always holds true.
Moreover, if in this case $\omega$ is the symplectic form on the nonsingular locus of $X$ induced by $f$ then the pullback of $\omega$ clearly extends to the original symplectic form on $\tilde X$.
Note that this is the usual definition for a resolution of singularities $f: \tilde X\to X$ of a symplectic variety $X$ to be symplectic.
\begin{Def}
Let $X$ be a scheme.
A (symplectic, $\Q$-factorial, ...) normal quasiprojective variety $\tilde X$ with at most terminal singularities together with a proper birational morphism $f: \tilde X\to X$ is called a (symplectic, $\Q$-factorial, ...) terminalisation (of $X$).
\end{Def}

\subsection{The Mukai vector}\label{DeltaMukaiSect}
Throughout this article $X$ will denote a projective K3 or abelian surface, which has Todd class $\td(X)=(1,0,2\epsilon)\in H^{2*}(X,\Z)$ with $\epsilon=1$ if $X$ is a K3 surface and $\epsilon=0$ if $X$ is an abelian surface.
The discriminant of a coherent sheaf $E$ on $X$ is
\begin{eqnarray}
\Delta(E):=2\rk E \c_2(E)-(\rk E-1)\c_1(E)^2 
= \c_1(E)^2 - 2\rk E \ch_2(E) \label{DeltaChern}
= v(E)^2 + 2\epsilon(\rk E)^2 \label{DeltaMukai}
\end{eqnarray}
and its Mukai vector is
$$v(E):=\ch(E)\sqrt{\td(X)}=(\rk E,\c_1(E),\chi(E)-\epsilon\rk E)\in \Hev:=\N_0\oplus\NS(X)\oplus\Z\,.$$
A vector $v$ of a lattice $\Lambda$ is primitive if there is no decomposition $v=mw$ with $2\le m\in\N$ and $w\in\Lambda$.
One has the even integral Mukai pairing
$$\langle (v_0,v_1,v_2),(v_0',v_1',v_2')\rangle := v_1.v_1'-v_0v_2'-v_0'v_2 $$
on $\Hev$.
We use the notation $v^2:=\langle v,v\rangle$.
By Mukai a simple sheaf on $X$ has always a Mukai vector $v$ with $v^2\ge -2$.

\subsection{General ample divisors}

The ample cone of $X$ carries a chamber structure for given Mukai vector $v=(v_0,v_1,v_2)\in\Hev$.
The definition depends on $v_0$. In the case of $v_0=1$ we agree that there is only one chamber, which equals the whole ample cone.

For $v_0>1$, we follow the definition in \cite{HL} section 4.C.
Let $\Num(X):=\Pic(X)/\equiv$, where $\equiv$ denotes numerical equivalence, and $\Delta:=v^2+2\epsilon v_0^2>0$ (this is the discriminant of a sheaf with Mukai vector $v$).
\begin{Def}
Let
$$W(v_0,\Delta):=\{ \xi^\perp \cap \Amp(X)_{\Q} \;|\; \xi\in\Num(X) \mwith -\frac {v_0^2}4 \Delta \le \xi^2 < 0 \}\,,$$
whose elements are called $v$-walls.
The connected components of the complement of the union of all $v$-walls are called $v$-chambers.
An ample divisor is called $v$-general if it is not contained in a $v$-wall.
\end{Def}
The set $W(v_0,\Delta)$ is locally finite in $\Amp(X)_{\Q}$ by \cite{HL} lemma 4.C.2.

For $v_0=0$, we follow the definition in \cite{Yo1} section 1.4.
\begin{Def}
Let $v_1\neq 0$ be effective. For every sheaf $F$ with $v(F)=v$ and every subsheaf $E\subset F$ we define $L:=\chi(E)\c_1(F)-\chi(F)\c_1(E)$, and for $L\ne 0$ we call
$$W_L:=L^\perp \cap\Amp(X)_\Q$$
the $v$-wall defined by $L$.
The connected components of the complement of the union of all $v$-walls are called $v$-chambers.
An ample divisor is called $v$-general if it is not contained in a $v$-wall.
\end{Def}
The number of nonempty $v$-walls is finite for a given Mukai vector $v=(0,v_1,v_2)$ with $v_1\neq 0$ effective by \cite{Yo1} section 1.4. (Yoshiokas additional assumption $v_1^2>0$ can be easily removed).

If $v_0=0=v_2$ then the notion of $H$-(semi)stability for a sheaf with Mukai vector $v$ is independent of the choice of $H$ and one cannot introduce the notion of a $v$-general ample divisor in this particular case.
However, we can move away from this case, as tensoring with the ample line bundle $H$ yields the isomorphism
$M_H(v)\cong M_H(v.\ch(H))\,.$
Thus one can assume without loss of generality that $v_2\neq 0$ when investigating the moduli spaces of one-dimensional semistable sheaves on a surface.

\section{A decomposition}

In this section, we exhibit the structure of the irreducible components of the moduli space $M_H(v)$. Although in general this moduli space is well-known to be irreducible, we could not extend this result to every case. On the other side, we have no reducible example at hand. Remark that, for an irreducible $M_H(v)$, we need these results for those moduli spaces containing no stable sheaves.

\begin{Prop}\label{MZerl}
Let $v\in\Hev$, $H$ an ample divisor and $M$ an irreducible component of $M_{H}(v)$.
Then there is a birational projective morphism $$g: \prod_{i=1}^m S^{n_i}M_i\to M$$
for a suitable decomposition $v=\sum_{i=1}^m n_iv_i$ with $n_i\in\N$ and $v_i\in\Hev$ for $1\le i\le m$ and a suitable choice of pairwise distinct irreducible components $M_i\subset \overline{M_H^s(v_i)}$ for $1\le i\le m$.
%
\end{Prop}
\Proof
Let $S$ be the at most countable set of finite tuples $(n_i,M_i)_i$ of pairs of a natural number $n_i\in\N$ and pairwise distinct connected components $M_i\subset \overline{M_H^s(v_i)}$ for some $v_i\in\Hev$ such that
there is an $H$-polystable sheaf $\left[\bigoplus_i \bigoplus_{j=1}^{n_i} F_{ij}\right]\in M$
with $F_{ij}\in M_i^s$ for all $1\le j\le n_i$ and all $i$.
For every such tuple $t=(n_i,M_i)_i\in S$ consider the morphism
$$g_t: \prod_i S^{n_i}M_i\to M_{H}(v),\ \left([F_{ij}]\right)\mapsto \left[\bigoplus_i \bigoplus_{j=1}^{n_i} F_{ij}\right]\,,$$
which is injective on the open subset ${ \prod_i S^{n_i}{M^s_i} }$.
The union of the images of the morphisms $g_t$ is $M_H(v)$. Since $t$ varies in a countable set, there exists a $\overline t$ such that $g_{\overline t}$ surjects on $M$.
Thus $g_{\overline t}$ is the desired birational morphism.
\qed\\

\noindent This result can be used to reduce the question of symplectic resolvability of any component of the moduli space of $H$-semistable sheaves to the question of symplectic resolvability of components of other moduli spaces of $H$-semistable sheaves containing stable sheaves:

\begin{Th}\label{MHss}
Let $v\in\Hev$, $H$ an ample divisor on $X$ and $M$ an irreducible component of $M_{H}(v)$.
Furthermore, let $$g: \prod_{i=1}^m S^{n_i}M_i\to M$$ be the projective birational morphism given by proposition \ref{MZerl}, where $M_i\subset \overline{M_H^s(v_i)}$ are pairwise distinct irreducible components.
Assume that for each $i$ there is a $\Q$-factorial projective symplectic terminalisation
$\tilde M_i\to M_i$.
For all $i$ set $$M_i^{(n_i)}:=\left\{\begin{array}{ll} \Hilb^{n_i}(\tilde M_i)&\mathrm{\ if\ } v_i^2=0,\\ S^{n_i} \tilde M_i&\mathrm{\ otherwise,}\end{array}\right.$$
and let $f$ be the concatenation of the projective birational morphisms
$$\tilde M:=\prod_{i=1}^m M_i^{(n_i)}\to\prod_{i=1}^m S^{n_i} \tilde M_i\to\prod_{i=1}^m S^{n_i} M_i\to M\,.$$
\begin{enumerate}
\item
If $\tilde M_i$ is nonsingular for all $i$ and $v_i^2\le 0$ whenever $n_i>1$
then $f$ is a projective symplectic resolution.

If $M'\to M$ is another projective symplectic resolution with $M'$ an irreducible symplectic manifold then it is deformation equivalent to $\tilde M_i$ for some $i$ or to a Hilbert scheme of points on a K3 surface.
\item 
If $\tilde M_j$ is singular or $v_j^2\ge 2$ and $n_j>1$ for some $j$
then $f$ is a singular $\Q$-factorial projective symplectic terminalisation.
\end{enumerate}
\end{Th}
\Proof
By \cite{Mu2}, since $M_i\neq\emptyset$, $M_i$ consists of one reduced point if $v_i^2<0$.
One easily verifies that $\tilde M$ is a projective symplectic variety with at most terminal singularities using \cite{Na4} corollary 1.
\begin{enumerate}
\item If $\tilde M_i$ is nonsingular for all $i$ and $v_i^2\le 0$ whenever $n_i>1$ then $\tilde M$ is nonsingular.

If $M'\to M$ is any other projective birational morphism with $M'$ an irreducible symplectic manifold then $\tilde M$ and $M'$ are deformation equivalent irreducible symplectic manifolds by a result of Huybrechts \cite{Huy}.
Furthermore, there is at most one $j$ with $v_j^2\ge 0$ and one has $n_j=1$ or $v_j^2=0$ for such a $j$.
In the second case $\tilde M=\Hilb^{n_j}(\tilde M_j)$, so $\tilde M_j$ must be a K3 surface.
\item If $\tilde M_j$ is singular or $v_j^2\ge 2$ and $n_j>1$ for some $j$
then $\tilde M$ is singular and $\Q$-factorial using the result of Bossi\`ere, Gabber and Serman that the direct product of two $\Q$-factorial varieties (over $\C$) is again $\Q$-factorial \cite{BGS}.\qed
\end{enumerate}

\section{Moduli spaces of one-dimensional sheaves}

The aim of this section is to extend the results of \cite{KLS} on moduli spaces for one-dimensional sheaves to all ample divisors.
First we connect the stability notion for varying ample divisors. 

\begin{Lemma}\label{D1Stab}
Let $v=(0,v_1,v_2)\in\Hev$ with $v_1\neq 0$, $v_2\neq 0$, $H$ an ample divisor in a $v$-chamber $K$ and $H'$ another ample divisor in the closure $\overline K$ of $K$ in the ample cone.
Then one has $H'$-stable $\Rightarrow$ $H$-stable $\Rightarrow$ $H$-semistable $\Rightarrow$ $H'$-semistable
for one-dimensional sheaves $F$ with $v(F)=v$.
\end{Lemma}
\Proof
Let $E\subset F$ be a nontrivial proper subsheaf of a sheaf $F$ with $v(F)=v$ and
$$f: \overline K\to \Q, h\mapsto (\chi(E)\c_1(F) - \chi(F)\c_1(E)).h\,. $$
Note that $f$ is well-defined on $\Q$-divisors already.
We consider the following two cases:
\begin{enumerate}
\item $F$ is $H$-semistable, i.e.\ $f(H)\le 0$, and we assume, $f(H') > 0$, or
\item $F$ is $H'$-stable, i.e.\ $f(H') < 0$, and we assume, $f(H)\ge 0$.
\end{enumerate}
Then there is a $\Q$-divisor $H_0$ on the connecting line of $H$ and $H'$ and $f(H_0) = 0$, and in particular, $H_0$ is $v$-general.
By definition of a $v$-general ample divisor, one thus has the contradiction $f\equiv 0$.
\qed\\

\noindent Now we can give the relations between the corresponding moduli spaces and state the consequences on the existence of symplectic resolutions. 

\begin{Th}\label{MHng0}
Let $v=(0,v_1,v_2)\in\Hev$ with $v_1\neq 0$, $v_2\neq 0$ and $v^2\ge 0$, and $H$ an ample divisor on $X$.
Assume that $M^s_H(v)$ is nonempty.
\begin{enumerate}
\item Let $v$ be primitive or $v^2=8$.
Then there is a projective symplectic resolution
$M\to \overline{M^s_H(v)}$.
If $H$ is not $v$-general then $M$ can be chosen to be a symplectic resolution of $M_A(v)$, where $A$ is a $v$-general ample divisor in a chamber touching $H$.
\item Let $v$ be not primitive and $v^2\neq 8$.
Then there is a singular locally factorial (and therefore $\Q$-factorial) projective symplectic terminalisation 
of $\overline{M^s_H(v)}\,.$
\end{enumerate}
\end{Th}
\Proof
For a $v$-general ample divisor $H$ this is due to Mukai, O'Grady, Kaledin, Lehn and Sorger, see e.g.\ \cite{KLS}, and for a not $v$-general ample divisor $H$ choose a $v$-general ample divisor $A$ in a chamber touching $H$ and use the projective birational morphism $f: M_A(v)\to M_H(v)$ induced by the universal properties of the moduli spaces thanks to the preceding lemma \ref{D1Stab}.
\qed

\section{Moduli spaces of two-dimensional sheaves}
For the rest of this article we restrict to sheaves of positive rank.
In the appendix we introduce the notion of (semi)stability with respect to a pair of two ample divisors $(H,A)$, which can be restated for two-dimensional sheaves on surfaces as follows:
A coherent sheaf $F$ is $(H,A)$-(semi)stable if it is $H$-semistable and if for any proper nontrivial subsheaf $E\subset F$ with reduced Hilbert polynomial $p_H(E)=p_H(F)$ one has $\mu_A(E)$ $(\ge)$ $\mu_A(F)$, i.e.\ stable corresponds to $>$ and semistable to $\ge$.

Let $H$ and $A$ be two ample divisors on $X$, $v=(v_0,v_1,v_2)\in\Hevp:=\N\oplus\NS(X)\oplus\Z$ on $X$ and
$P$ the Hilbert polynomial of any sheaf with Mukai vector $v$.
We define $M_{H,A}(v)$ to be the fibre of the morphism
$M_{H,A}(P) \stackrel{\det}{\to} \Pic(X) \stackrel{\c_1}{\to} \NS(X)$
over $v_1\in\NS(X)$.
This moduli space parametrises $(H,A)$-polystable sheaves with Mukai vector $v$.
We denote the open subset of $(H,A)$-stable sheaves by $M^s_{H,A}(v)$.
The local description of the moduli space $M_{H,A}(v)$ is analogous to $M_H(v)$.
In particular, $M^s_{H,A}(v)$ is nonsingular, each connected component has even dimension $2+v^2$, and Mukai's construction of the symplectic form carries over to it.

The following proposition points out the importance of $v$-general ample divisors. The proportionality in the first part is crucial for the local analysis of the moduli spaces.

\begin{Prop}\label{propp}
Let $A$ be $v$-general and $F$ an $(H,A)$-semistable sheaf with Mukai vector $v$.
\begin{enumerate}
\item
If $E\subset F$ is a nontrivial proper subsheaf with $p_{H,A}(E)=p_{H,A}(F)$ then
$\frac {\c_1(E)}{\rk E}=\frac {\c_1(F)}{\rk F}\,.$
\item
If $v$ is primitive then $F$ is $(H,A)$-stable.
\end{enumerate}
\end{Prop}
\Proof
\begin{enumerate}
\item
$E$ is saturated, so $0<\rk E<\rk F$.
By lemma \ref{StabChar2d} one has
$$\left( \mu_H(E), \frac{\chi(E)}{\rk E}, -\mu_A(E) \right)=\left( \mu_H(F), \frac{\chi(F)}{\rk F}, -\mu_A(F) \right) \,.$$
Assume that $\xi:=\rk F\;\c_1(E)-\rk E\;\c_1(F)\neq 0$.
As $F$ is $\mu_H$-semistable, by \cite{HL} theorem 4.C.3 $\xi$ defines the $v$-wall $\xi^\perp\cap\Amp(X)_\Q$. This wall contains $A$, a contradiction to $A$ being $v$-general.
\item
Assume $E\subset F$ is a nontrivial saturated proper subsheaf with $p_{H,A}(E)=p_{H,A}(F)$, so
$\rk E(\c_1(F),\chi(F))=\rk F(\c_1(E),\chi(E))$
by item 1. Thus
$$\frac {\rk E}{\gcd(\rk E,\rk F)}(\c_1(F),\chi(F))=\frac {\rk F}{\gcd(\rk E,\rk F)}(\c_1(E),\chi(E))\,.$$
Clearly $1\le \gcd(\rk E,\rk F)\le \rk E<\rk F\,,$ hence  $\frac {\rk F}{\gcd(\rk E,\rk F)}>1\,$ 
and this integer divides $(\c_1(F),\chi(F))$ and $\rk F$, which is a contradiction to the primitiveness of $v$.
\end{enumerate}
The following result holds even for arbitrary ample divisors $A$:
\begin{samepage}
\begin{Th}\label{StabKompHA} 
If $M\subseteq M_{H,A}(v)$ is a connected component with $M\subseteq M^s_{H,A}(v)$ then one already has $M=M_{H,A}(v)$.
In particular, if $M^s_{H,A}(v)=M_{H,A}(v)$ then $M_{H,A}(v)$ is irreducible.
\end{Th}
\Proof
The proof of \cite{KLS} theorem 4.1 carries over literally.
\qed\\
\end{samepage}

\noindent We come to our main result on the the moduli space of $(H,A)$-semistable sheaves.

\begin{Th}\label{MAHgen}
Let $v=(v_0,v_1,v_2)\in\Hevp$ be primitive, $m\in\N$, $H$ an ample divisor on $X$ and $A$ an $mv$-general ample divisor on $X$.
Assume that $M^s_{H,A}(mv)$ is nonempty.
\begin{enumerate}
\item If $v^2=0$ then $M_{H,A}(mv)=M^s_{H,A}(mv)$, and $M_{H,A}(mv)$ is a projective symplectic nonsingular surface.
\item Let $v^2\ge 2$ and $M^s_{H,A}(v)$ be nonempty.
Then $M_{H,A}(mv)$ is a projective symplectic variety of dimension $2+m^2v^2$.
\begin{enumerate}
\item If $m=1$ then $M_{H,A}(v)=M^s_{H,A}(v)$, and $M_{H,A}(v)$ is nonsingular.
\item If $m\ge 2$ then the singular locus of $M_{H,A}(mv)$ is nonempty and equals the strictly semistable locus.
\begin{enumerate}
\item If $m=2$ and $v^2=2$ then the singular locus has codimension $2$ and $M_{H,A}(mv)$ admits a symplectic resolution.
\item If $m=2$ and $v^2>2$ or $m>2$ then $M_{H,A}(mv)$ is locally factorial, the singular locus has codimension at least $4$ and the singularities are terminal. There is no open neighbourhood of a singular point that admits a symplectic resolution.
\end{enumerate}
\end{enumerate}
\item Let $v^2\ge 2$ but now $M^s_{H,A}(v)$ be empty. Then $M_{H,A}(v)$ is empty as well, i.e.\ $m>1$ by assumption.
If $m=2$ or $3$ then $M_{H,A}(mv)=M^s_{H,A}(mv)$, and $M_{H,A}(mv)$ is a nonsingular projective symplectic variety of dimension $2+m^2v^2$.
\end{enumerate}
\end{Th}
The third part of this theorem is only included as our existence results for stable sheaves presented in the following section do not fully extend the existence results one has for $H$-stable sheaves. We don't have an example at hand for this case.\\

\ProofOf \textit{the theorem}.
If $m=1$ then $M^s_{H,A}(v)=M_{H,A}(v)$ by proposition \ref{propp}.
If $M^s_{H,A}(mv)=M_{H,A}(mv)$ then $M_{H,A}(mv)$ is irreducible by theorem \ref{StabKompHA}. In particular, $M_{H,A}(v)$ is a projective symplectic nonsingular variety.

\begin{enumerate}
\item Let $v^2=0$.
As $M^s_{H,A}(mv)$ is nonempty, $M^s_{H,A}(av)$ is empty for all $a\neq m$ analogously to the arguments in $\cite{Mu2}\ \S 3$.
Thus $M_{H,A}(mv)=M^s_{H,A}(mv)$.
\item Let $v^2\ge2$ and $M^s_{H,A}(v)$ be nonempty.
The case of $m=1$ is clear by the above statements, so let $m\ge 2$.
The results herein are straightforward generalisations of \cite{KLS} and carry over literally as
they are based on a local analysis using proposition \ref{propp},
and the local description of $M_{H,A}(mv)$ is analogous to the one of $M_A(mv)$.
\item Let $v^2\ge2$ and $M^s_{H,A}(v)$ be empty. If $m=2$ or $3$ then $M_{H,A}(mv)=M^s_{H,A}(mv)$, and the claim follows from the statements at the beginning of the proof.
\qed
\end{enumerate}
The moduli spaces $M_{H,A}(v)$ and their symplectic resolutions are good candidates for $\Q$-factorial projective symplectic terminalisations of $M_{H}(v)$,
and we can reduce our question on $M_H(v)$ to the investigation of $M_{H,A}(v)$:

\begin{Prop}\label{MHngPTerm}
Let $v\in\Hevp$, $H$ a not $v$-general and $A$ a $v$-general ample divisor.
Assume that $M^s_H(v)$ is nonempty and that there is a $\Q$-factorial projective symplectic terminalisation
$M\to M_{H,A}(v)\,.$
Then there is a $\Q$-factorial projective symplectic terminalisation $f: M\to\overline{M^s_H(v)}\,.$
\end{Prop}
\Proof
Concatenate the terminalisation $M\to M_{H,A}(v)$ with the morphism $M_{H,A}(v) \to \overline{M^s_H(v)}$ induced by the universal properties of the moduli spaces.\qed\\

\noindent Since we were unable to extend the irreducibility result to all moduli spaces of $(H,A)$-semistable sheaves, as a substitute we consider components containing no $(H,A)$-stable sheaves. These considerations are also relevant for the irreducible moduli spaces containing no stable sheaves at all, e.g.\ in the case of an isotropic Mukai vector.

\begin{Cor}\label{MHssP}
Let $v=(v_0,v_1,v_2)\in\Hevp$, $H$ a not $v$-general ample divisor on $X$ and $M$ an irreducible component of $M_{H}(v)$ containing no $H$-stable sheaves.

Assume that for all $w=(w_0,w_1,w_2)\in\Hevp$ with 
$1<w_0<v_0$ and such that $H$ is not $w$-general, $\frac{w_1.H}{w_0}=\frac{v_1.H}{v_0}$ and
$\frac{w_2}{w_0}=\frac{v_2}{v_0}$
there is a $\Q$-factorial projective symplectic terminalisation of $M_{H,A_w}(w)$
for a suitable $w$-general ample divisor $A_w$.
Then there is a $\Q$-factorial projective symplectic terminalisation $\tilde M\to M$.

If $\tilde M$ can be chosen to be an irreducible symplectic manifold then it is deformation equivalent
to some symplectic resolution of some $M_{H,A}(w)$, where $w=(w_0,w_1,w_2)\in\Hevp$ has the above properties and $A$ is a $w$-general ample divisor,
to a symplectic resolution of some $M_H(w)$, where $1\le w_0<v_0$, $H$ is $w$-general, $\frac{w_1.H}{w_0}=\frac{v_1.H}{v_0}$ and $\frac{w_2}{w_0}=\frac{v_2}{v_0}$,
or to a Hilbert scheme of points on a K3 surface.
\end{Cor}
\Proof
Consider the decomposition $v=\sum_{i=1}^m n_iv^{(i)}$ given by proposition \ref{MZerl}.
As the Mukai vectors $v^{(i)}$ belong to $H$-stable direct summands of a strictly $H$-polystable sheaf with Mukai vector $v$ one has 
$1<v^{(i)}_0<v_0$, $\frac{v^{(i)}_1.H}{v^{(i)}_0}=\frac{v_1.H}{v_0} \mand \frac{v^{(i)}_2}{v^{(i)}_0}=\frac{v_2}{v_0}$ for all $i$.
If $H$ is $v^{(i)}$-general for some $i$ then $M_H(v^{(i)})$ is a symplectic variety that admits a $\Q$-factorial projective symplectic terminalisation by \cite{KLS}.
Thus theorem \ref{MHss} yields the claim.
\qed

\section{Existence of stable sheaves}
Given a primitive $v=(v_0,v_1,v_2)\in\Hevp$, $m\in\N$, $H$ an ample divisor on $X$ and $A$ an $mv$-general ample divisor on $X$, theorem \ref{MAHgen} extends well-known results on $M_A(mv)$ to $M_{H,A}(mv)$ assuming the existence of $(H,A)$-stable sheaves with Mukai vector $v$.
In this section we deduce an existence result for these sheaves stated in theorem \ref{MAHBstab}. First, we need some preliminaries.

\begin{Lemma}\label{FiltDeltaEx}
Let $0=F_0\subset F_1\subset ... \subset F_n = F$ be a filtration of a coherent sheaf $F$ on $X$ with positive rank such that the graded objects $gr_i:=F_i/F_{i-1}$ have positive rank for $i=1,...,n$.
Then
$$\sum_{i=1}^n \frac {\Delta(gr_i)} {\rk gr_i} - \frac {\Delta(F)} {\rk F}=
\sum_{i<j} \frac {\rk gr_i\ \rk gr_j} {\rk F} \left(\frac {\c_1(gr_i)} {\rk gr_i}-\frac {\c_1(gr_j)} {\rk gr_j}\right)^2\,.$$
\end{Lemma}
\Proof See the proof of \cite{HL} corollary 7.3.2.
\qed

\begin{Cor}\label{FiltDelta}
If all $gr_i$ have the same slope with respect to $H$ then one has
$$\sum_{i=1}^n \frac {\Delta(gr_i)} {\rk gr_i} \le \frac {\Delta(F)} {\rk F}\,.$$
Moreover, if $\frac {\c_1(gr_i)} {\rk gr_i}\neq\frac {\c_1(gr_j)} {\rk gr_j}$ for all $1\le i<j\le n$ then one even has
\begin{eqnarray*}
\sum_{i=1}^n \frac {\Delta(gr_i)} {\rk gr_i}
&\le& \frac {\Delta(F)} {\rk F}-\sum_{i<j} \frac {2\rk gr_i \rk gr_j}{\rk F\; \lcm(\rk gr_i,\rk gr_j)^2}
\le \frac {\Delta(F)} {\rk F}-\sum_{i<j} \frac {2}{\rk F \rk gr_i \rk gr_j}\,,
\end{eqnarray*}
where $\lcm$ denotes the least common multiple.
\end{Cor}
\Proof
By assumption one has
$\left(\frac {\c_1(gr_i)} {\rk gr_i}-\frac {\c_1(gr_j)} {\rk gr_j}\right).H=0\,,$
hence $\left(\frac {\c_1(gr_i)} {\rk gr_i}-\frac {\c_1(gr_j)} {\rk gr_j}\right)^2\le 0$ for all $i,j$ by the Hodge index theorem.
The intersection pairing is nondegenerate and even, hence
\begin{eqnarray*}
\left(\frac {\c_1(gr_i)} {\rk gr_i}-\frac {\c_1(gr_j)} {\rk gr_j}\right)^2
&=&\frac 1{\lcm(\rk gr_i,\rk gr_j)^2}\left(\lcm(\rk gr_i,\rk gr_j)\left(\frac {\c_1(gr_i)} {\rk gr_i}-\frac {\c_1(gr_j)} {\rk gr_j}\right)\right)^2\\
&\le& -\frac {2}{\lcm(\rk gr_i,\rk gr_j)^2}\,.
\end{eqnarray*}
\qed

\begin{Lemma}\label{ChiIneq}
Let $H$ be an ample divisor on $X$, $2\le n\in\N$ and
$0=F_0\subset F_1\subset ... \subset F_n = F$ a filtration of a coherent sheaf $F$ on $X$ with positive rank $r$ such that all $gr_i=F_i/F_{i-1}$ have positive rank $r_i$, are $\mu_H$-semistable, have the same slope with respect to $H$ and
$\frac {\c_1(gr_i)} {\rk gr_i}\neq\frac {\c_1(gr_j)} {\rk gr_j}$ for all $1\le i<j\le n$.
Then
$$\sum_{i<j} \chi(gr_i,gr_j) \le - \frac{\Delta(F)}{2r}(n-1) + \epsilon r^2 -\epsilon\sum_{i=1}^n r_i^2-\frac{r-n+1}{r}\sum_{i<j} \frac 1{r_ir_j}\,.$$
\end{Lemma}
\Proof
One has $r_i\le r-n+1$ and the Bogomolov inequality $\Delta(gr_i)\ge 0$ for all $i$.
Using equation (\ref{DeltaMukai}) of section \ref{DeltaMukaiSect} we can calculate
\begin{eqnarray*}
\sum_{i<j} \chi(gr_i,gr_j)&=&-\sum_{i<j} \langle v(gr_i),v(gr_j)\rangle\\
&=&\frac 1 2\left(-\sum_{i,j=1}^n \langle v(gr_i),v(gr_j)\rangle + \sum_{i=1}^n v(gr_i)^2\right)\\
&=&\frac 1 2\left(\sum_{i=1}^n v(gr_i)^2 - v(F)^2 \right)\\
&\stackrel{\mathrm{(\ref{DeltaMukai})}}=& \frac 1 2\left(\sum_{i=1}^n \left(\Delta(gr_i)- 2\epsilon r_i^2\right) - \Delta(F) + 2\epsilon r^2\right)\\
&\le& \sum_{i=1}^n \frac{r-n+1}{2r_i}\Delta(gr_i)- \epsilon\sum_{i=1}^n r_i^2 - \frac {\Delta(F)} 2 + \epsilon r^2\\
&\stackrel{\mathrm{cor.\ \ref{FiltDelta}}}\le& \frac{r-n+1}2\left( \frac{\Delta(F)}{r} -\sum_{i<j} \frac {2}{rr_ir_j}\right)-\epsilon\sum_{i=1}^n r_i^2 - \frac {\Delta(F)} 2 + \epsilon r^2\\
&=& - \frac{\Delta(F)}{2r}(n-1) + \epsilon r^2 -\epsilon\sum_{i=1}^n r_i^2-\frac{r-n+1}{r}\sum_{i<j} \frac 1{r_ir_j}\,.
\end{eqnarray*}
\qed\\
\pagebreak[2]

\noindent Some parts of the proof of the following proposition are based on an idea we learned from an unpublished note of Christoph Sorger.

\begin{Prop}\label{ABstnHst}
Let $v=(v_0,v_1,v_2)\in\Hevp$ with $v_0\ge 2$ and
$H$, $A$ and $B$ three ample divisors on $X$.
Assume that $M^s_{A,B}(v)$ is nonempty and contains no $H$-semistable sheaves,
and let $R^s\to M^s_{A,B}(v)$ be the geometric quotient of the construction of the moduli space $M_{A,B}(v)$ and $F\in\Coh(R^s\times X)$ the universal quotient family.
Then there is an open dense subset $S\subset R^s$ and a subsheaf $F'\subset F|_S$ such that for all $s\in S$ one has
\begin{enumerate}
\item an exact sequence $0\to F'_s\to F_s\to F''_s\to 0$ on the fibre over $s$ with
\item $p_H(F'_s)>p_H(F_s)>p_H(F''_s)$,
\item $\hom(F'_s,F''_s)=0$ and
\item $1-\chi(F'_s,F''_s) = \ext^2_-(F_s,F_s)\le \emor(F'_s)+\emor(F''_s)$,
\end{enumerate}
where we calculate $\ext^2_-(F_s,F_s)$ with respect to the filtration $F'_s\subset F_s$
(for a definition see \cite{HL} section 2.A)
%
,
and if $v_0=2$ then additionally
\begin{enumerate}
\item $\ext^2_-(F_s,F_s)=2$,
\item $-\chi(F'_s,F''_s)=\ext^1(F'_s,F''_s)=1$ and
\item $0\ge \left(\c_1(F'_s) - \c_1(F''_s) \right)^2 = v^2-4 = v(F'_s)^2 + v(F''_s)^2 -2 \,.$
\end{enumerate}
\end{Prop}
\Proof
Stable sheaves are simple, hence one has
\begin{eqnarray}
\ext^2(F_s,F_s)=\hom(F_s,F_s)=1\label{MAHB1}
\end{eqnarray}
for all $s\in S$.
By the same arguments as in the proof of \cite{HL} theorem 10.2.1 
$R^s$ is nonsingular and the Kodaira-Spencer map $\kappa$ is given by the concatenation of the two maps
$$T_{s}R^s\stackrel{}\longrightarrow T_{[F_s]} M^s_A(v)\stackrel{\cong}\longrightarrow\Ext^1(F_s,F_s)\,.$$
Furthermore, the first map is surjective, hence $\kappa$ is surjective as well.

In the following every notion is understood to be with respect to the ample divisor $H$ whenever not explicitly stated differently.
By \cite{HL} theorem 2.3.2 there is a relative Harder-Narasimhan filtration $F_\bullet$ and an open dense subscheme $S\subset R^s$ such that the restriction of the filtration to a fibre over $s\in S$ is a Harder-Narasimhan filtration of $F_s$.
As the open subset of $R^s$ containing $H$-semistable sheaves is empty the filtration is nontrivial.
We only take the first step $F':=F_{\ell-1}|_S\subset F_\ell|_S=F|_S$ of the filtration restricted to $S$, which gives us an exact sequence
\begin{eqnarray}
0\to F'_s\to F_s\to F''_s\to 0\label{Fseq}
\end{eqnarray}
on the fibres over $s\in S$ with
\begin{eqnarray}
p_H(F'_s)>p_H(F_s)>p_H(F''_s)\,.\label{MAHB11}
\end{eqnarray}
By the proof of \cite{HL} theorem 2.3.2 one has
\begin{eqnarray}
\hom(F'_s,F''_s)&=&0\label{MAHB2}
\end{eqnarray}
and
$\pi: \Quot_{S\times X/S}(F,P_-)\to S$
is an isomorphism, where $P_-:=P_H(F''_s)$ (this is independent of $s$).
Let $s\in S$ be a closed point and $x$ be the unique point with $s=\pi(x)$,
which corresponds to the exact sequence \ref{Fseq}.
By \cite{HL} theorem 2.2.7 the kernel of the obstruction map $\mfo: T_sS\to \Ext^1(F'_s,F''_s)$
is $\ker\mfo = \Im\;T_x\pi = \dim T_sS\,,$
hence $\mfo$ is the zero map.
As $\mfo$ is given by
$$\mfo: T_sS \stackrel{\kappa}\longrightarrow \Ext^1(F_s,F_s) \stackrel {c}\longrightarrow \Ext^1(F'_s,F''_s)$$
and $\kappa$ is surjective as explained above one has $c=0$ as well.
For the short filtration $0\subset F'_s\subset F_s$ there is a long exact sequence
$$...\to \Ext^i_-(F_s,F_s) \to \Ext^i(F_s,F_s)\to \Ext^i(F'_s,F''_s) \to  \Ext^{i+1}_-(F_s,F_s)\to ... \,,$$
which decomposes to the exact sequence
$$ 0\to \Ext^1(F'_s,F''_s) \to \Ext^2_-(F_s,F_s) \to \Ext^2(F_s,F_s)\to \Ext^2(F'_s,F''_s) \to 0$$
as $c=0$. Thus
\begin{eqnarray*}
0&=& -\ext^1(F'_s,F''_s) + \ext^2_-(F_s,F_s) - \ext^2(F_s,F_s)+ \ext^2(F'_s,F''_s)\\
&\stackrel{\mathrm{(\ref{MAHB1}),\ (\ref{MAHB2})}}=& \chi(F'_s,F''_s) + \ext^2_-(F_s,F_s) -1\,.
\end{eqnarray*}
By \cite{HL} theorem 2.A.4
there is a spectral sequence
\begin{eqnarray*}
\Ext^{p+q}_-(F_s,F_s) &\Leftarrow& E^{pq}_1=\left\{\begin{array}{ll} 0& p<0\\ \prod_i\Ext^{p+q}(gr_i F_s,gr_{i-p} F_s) &p\ge 0\end{array}\right.\,.
\end{eqnarray*}
Hence
\begin{eqnarray}
\ext^2_-(F_s,F_s)&\le& \sum_{i\ge j} \ext^2(gr_i F_s,gr_j F_s)\nonumber\\ 
&=& \sum_{i\ge j} \hom(gr_j F_s,gr_i F_s)\nonumber\\
&=& \emor(F'_s)+\emor(F''_s)+\hom(F'_s,F''_s)\nonumber\\
&\stackrel{\mathrm{(\ref{MAHB2})}}=& \emor(F'_s)+\emor(F''_s)\label{MAHB9}
\end{eqnarray}
Assume $v_0=2$.
Then $F'_s$ and $F''_s$ are line bundles and therefore $(A,B)$-stable and simple.
As $F_s$ is $(A,B)$-stable, one has
\begin{eqnarray}
p_{A,B}(F'_s) <_0 p_{A,B}(F_s) <_0 p_{A,B}(F''_s)\label{MAHB12}
\end{eqnarray}
and in particular,
$0=\hom(F''_s,F'_s)=\ext^2(F'_s,F''_s)$
analogous to \cite{HL} proposition 1.2.7.
The sequence $0\to F'_s\to F_s\to F''_s\to 0$ gives a nontrivial element in $\Ext^1(F'_s,F''_s)$, thus 
$\ext^1(F'_s,F''_s)\ge 1$.
Altogether one has
$$1\le \ext^1(F'_s,F''_s) = -\chi(F'_s,F''_s) = \ext^2_-(F_s,F_s)-1 \stackrel{\mathrm{(\ref{MAHB9})}}\le 1\,,$$
which gives us equality everywhere. In particular,
\begin{eqnarray*}
-2&=&2\chi(F'_s,F''_s)\\
&=&-2\langle v(F'_s),v(F''_s)\rangle\\
&=& v(F'_s)^2 + v(F''_s)^2 - v^2 \\
&\stackrel{\mathrm{(\ref{DeltaMukai})}}=& \Delta(F'_s)+\Delta(F''_s) - 4\epsilon- v^2\\
&\stackrel{\mathrm{l.\ \ref{FiltDeltaEx}}}=& \frac{\Delta(F_s)}2 + \frac 1 2 \left(\c_1(F'_s) - \c_1(F''_s) \right)^2  - 4\epsilon- v^2\\
&\stackrel{\mathrm{(\ref{DeltaMukai})}}=&\frac 1 2 \left(\c_1(F'_s) - \c_1(F''_s) \right)^2  - \frac {v^2}2\,.
\end{eqnarray*}
By the inequalities \ref{MAHB11} and \ref{MAHB12} one has
\begin{eqnarray*}
\mu_H(F'_s)&\ge&\mu_H(F_s)\mand\\
\mu_A(F'_s)&\le&\mu_A(F_s)\,,
\end{eqnarray*}
hence there is an ample $\Q$-divisor $H'\in [H,A]$ such that
$\mu_{H'}(F'_s)=\mu_{H'}(F_s)=\mu_{H'}(F''_s)\,,$
and the Hodge index theorem (see e.g.\ \cite{HL} theorem V.1.9.) yields
$$0\ge \left(\c_1(F'_s) - \c_1(F''_s) \right)^2 = v^2-4\,.$$
\qed\\

\noindent We are now ready for the main result of this section:

\begin{Th}\label{MAHBstab}
Let $v=(v_0,v_1,v_2)\in\Hevp$ with $v_0\ge 2$ and
\begin{eqnarray*}
v^2 > 2\left(  v_0^3-(2-\epsilon)v_0^2+v_0(1-2\epsilon)-(v_0-1)\left\lfloor \frac {v_0^2}4\right\rfloor^{-1}\right)=:\varphi(v_0)\,,
\end{eqnarray*}
$H$ a not $v$-general ample divisor and $A$ a $v$-general ample divisor in a chamber touching $H$.
Then there is an $A$-stable and $H$-semistable sheaf $E$ with Mukai vector $v$.

Moreover, if $B$ is another $v$-general ample divisor such that $H\in[A,B]$ is the unique not $v$-general ample divisor then $E$ is even $(H,B$)-stable.
\end{Th}
\Proof
One can easily verify that by the assumptions on $v$ one has $v^2\ge 2$, and therefore $M^s_A(v)\neq\emptyset$ by \cite{KLS}.
Set $B=A$ and assume that $M^s_A(v)$ contains no $H$-semistable sheaves.
Then by proposition \ref{ABstnHst} there is an $[F]\in M^s_A(v)$ with an exact sequence 
$$0\to F'\to F\to F''\to 0$$ such that
\begin{eqnarray}
&&p_H(F')>p_H(F)>p_H(F'')\,,\label{MAHB10}\\
&&\hom(F',F'')=0 \mand \\
&&1-\chi(F',F'') = \ext^2_-(F,F)\le \emor(F')+\emor(F'')\,. \label{MAHB5}
\end{eqnarray}
As $F$ is in particular $\mu_A$-semistable, it is also $\mu_H$-semistable.
Together with inequality (\ref{MAHB10}) one has
\begin{eqnarray}
&&\mu_H(F')=\mu_H(F)=\mu_H(F'')\label{MAHB6} \mand\\
&&\frac{\chi(F')}{\rk F'}>\frac{\chi(F)}{\rk F}>\frac{\chi(F'')}{\rk F''},\label{MAHB3}
\end{eqnarray}
and $F'$ and $F''$ are $\mu_H$-semistable.
Thus by \cite{O3} lemma 1.7 one has
\begin{eqnarray}
\emor(F')\le (\rk F')^2 &\mand&\emor(F'')\le (\rk F'')^2 \,.\label{MAHB7}
\end{eqnarray}
Moreover, the $A$-stability of $F$ ensures
$$p_A(F')<p_A(F)<p_A(F'')\,,$$
and because of inequality (\ref{MAHB3}) even
$$\mu_A(F')<\mu_A(F)<\mu_A(F'')\,,$$
and therefore
\begin{eqnarray}
\frac{\c_1(F')}{\rk F'}\neq\frac{\c_1(F'')}{\rk F''}\,.\label{MAHB4}
\end{eqnarray}
We will need the easily verified inequalities
\begin{eqnarray}
(\rk F')^2+(\rk F'')^2&\le& v_0^2-2v_0+2 \label{MAHB20} \mand  \\ 
\rk F'\; \rk F'' &\le& \left\lfloor \frac {v_0^2}4\right\rfloor \label{MAHB21} \,.
\end{eqnarray}
Altogether one has
\begin{eqnarray*} 
1&\stackrel{\mathrm{(\ref{MAHB5})}}\le & \chi(F',F'') + \emor(F')+\emor(F'')\\
&\stackrel{\mathrm{l.\ \ref{ChiIneq},\ (\ref{MAHB7})}}\le& 
- \frac{\Delta(F)}{2v_0} + \epsilon v_0^2 + (1-\epsilon)( (\rk F')^2+(\rk F'')^2) -\frac{v_0-1}{v_0\;\rk F'\; \rk F''}\\
&\stackrel{\mathrm{(\ref{DeltaMukai})}}\le&
- \frac{v^2}{2v_0} + \epsilon (v_0^2-v_0) + (1-\epsilon)( (\rk F')^2+(\rk F'')^2) -\frac{v_0-1}{v_0\;\rk F'\; \rk F''}\\
&\stackrel{(\ref{MAHB20}),(\ref{MAHB21})}{\le}&
- \frac{v^2}{2v_0} + \epsilon (v_0^2-v_0) + (1-\epsilon)( v_0^2-2v_0+2) -\frac{v_0-1}{v_0}\left\lfloor \frac {v_0^2}4\right\rfloor^{-1}\\
&=& -\frac{v^2}{2v_0}+ v_0^2-(2-\epsilon)v_0+2(1-\epsilon)-\frac{v_0-1}{v_0}\left\lfloor \frac {v_0^2}4\right\rfloor^{-1}\,,
\end{eqnarray*}
in contradiction to the assumption of the theorem.

Thus there is an $A$-stable and $H$-semistable sheaf $E$ with Mukai vector $v$.
Let $B$ be another $v$-general ample divisor such that $H\in[A,B]$ is the unique not $v$-general ample divisor.
Let $E'\subset E$ be a proper nontrivial subsheaf with $p_H(E')=p_H(E)$. In particular, 
$\mu_H(E')=\mu_H(E)$, hence the linear map
$$[A,B]\to\R,\quad h\mapsto \mu_h(E')-\mu_h(E)$$
is either zero everywhere or changes the sign when passing through $H$.
In the first case one has $\mu_A(E')=\mu_A(E)$ and therefore $\frac{\chi(E')}{\rk E'}<\frac{\chi(E)}{\rk E}$ by the $A$-stability of $E$.
In the second case one has $\mu_A(E')<\mu_A(E)$ and therefore $\mu_B(E')>\mu_B(E)$.
By the characterisation in lemma \ref{StabChar2d} one has that $p_{H,B}(E') <_0 p_{H,B}(E)$ in both cases.
\qed

\begin{Rem}
Since $B$-stability and $(H,A)$-stability are open conditions, the existence of a sheaf which is both $B$-stable and $(H,A)$-stable induces an isomorphism between an open subset of $M_B(v)$ and an open subset of $M_{H,A}(v)$.
\end{Rem}

\begin{Cor}\label{DefClass}
Let $X$ be a K3 surface, $v=(v_0,v_1,v_2)\in\Hevp$ a primitive Mukai vector, $H$ an ample divisor and $A$ a $v$-general ample divisor.
Assume that $v^2>\varphi(v_0)$.
Then $M_{H,A}(v)$ is deformation equivalent to $\Hilb^{\frac{v^2}2+1}(X)$.
\end{Cor}
\Proof
$M_{H,A}(v)$ is irreducible by theorem \ref{MAHgen}. Theorem \ref{MAHBstab} together with the above remark yields birationality of $M_{H,A}(v)$ and $M_B(v)$,
and the deformation equivalence follows by Huybrechts \cite{Huy} using the results given in the introduction of \cite{KLS}.
\qed\\

\noindent Let us evaluate $\varphi$ for small values:
\begin{center}
\begin{tabular}{|c|c|c|c|c|c|c|c|c|c|c|}
\hline
$n$&2&3&4&5&6\\
\hline
$\varphi(n)$ for $\epsilon$=0&2&22&70.5&158.7&298.9\\
\hline
$\varphi(n)$ for $\epsilon$=1& 2&28&86.5&188.7&346.9\\
\hline
\end{tabular}
\end{center}

\noindent 
In particular, the only interesting exceptional case for rank two might occur for $v^2=2$.
To realise this case one needs a K3 surface that holds a divisor $D$ with $D^2=-2$, $D.H=0$, $D.A<0$ and such that
$D+v_1$ is divisible by 2 in $\NS(X)$. This cannot occur as if $D$ is a divisor with $D^2=-2$ then $D$ or $-D$ is effective, hence $D.H\neq 0$ for any ample divisor $H$.

\vspace{1cm}
\begin{center}
\Large \bf Appendix
\end{center}
In this appendix we introduce the notion of (semi)stability of sheaves on a projective scheme with respect to a pair of two ample divisors $(H,A)$, give some properties and construct the corresponding moduli space.

\section{Semistable sheaves}\label{Semistable}
We assume familiarity with the material presented in \cite{HL} and use the notation therein.
Let $X$ be a projective scheme over a field $k$, $H$ and $A$ two ample line bundles on $X$ and $E$ a nontrivial coherent sheaf on $X$.
We write $E(mH+nA):=E\otimes H^{\otimes m} \otimes A^{\otimes n}\,,$
and we denote the Hilbert polynomial of $E$ with respect to $H$ by $P_H(E)$, its reduced Hilbert polynomial by $p_H(E)$ and its multiplicity by $\alpha^H_{\dim E}(E)$.
We define
$$
P_{H,A}(E)(m,n):=\chi(E(mH+nA)) \mand
p_{H,A}(E):=\frac{P_{H,A}(E)}{\alpha^H_{\dim E}(E)} \,.
$$
These are polynomials in $m$ and $n$ with degree $d:=\dim E$ in $n$ and $m$ and total degree $d$, and one has $P_{H,A}(E)(\bullet,0)=P_{H}(E)$ and $p_{H,A}(E)(\bullet,0)=p_{H}(E)$.

There is a natural ordering of polynomials in one variable given by the lexicographic ordering of their coefficients.
This generalises to polynomials of two variables by the identification $\Q[m,n]=(\Q[m])[n]$, i.e.\ we consider the elements as polynomials in $n$ and use the ordering of $\Q[m]$ for comparing coefficients.

We introduce another ordering on $\Q[m,n]$ by defining
$$f\le_0 g\quad:\Leftrightarrow\quad (f(\bullet,0),-f)\le(g(\bullet,0),-g)$$ for $f,g\in\Q[m,n]$, where on the right hand side we use lexicographic ordering on the product $\Q[m]\times\Q[m,n]$, i.e.\ $f\le_0 g$ if and only if $f(\bullet,0)<g(\bullet,0)$ or $f(\bullet,0)=g(\bullet,0)$ and $f\ge g$.
Clearly one has $f=_0g$ if and only if $f=g$.

In order to avoid case differentiation for stable and semistable sheaves we follow the notation 1.2.5 in \cite{HL} using bracketed inequality signs, e.g.\ an inequality with $\lel$ for (semi)stable sheaves means that one has $\le$ for semistable sheaves and $<$ for stable sheaves.
\begin{Def}\label{HAhst}
A coherent sheaf $E$ of dimension $d$ is $(H,A)$-(semi)stable if it is pure and if for any proper nontrivial subsheaf $F\subset E$ one has $p_{H,A}(F)\lelo p_{H,A}(E)\,.$
If $E$ is strictly $(H,A)$-semistable, i.e.\ $(H,A)$-semistable but not $(H,A)$-stable, then there is always a proper nontrivial subsheaf $F\subset E$ with $p_{H,A}(F)=p_{H,A}(E)$, which is then called an $(H,A)$-destabilising subsheaf.
\end{Def}
This definition is independent of the choice of the two ample line bundles in $\Q\cdot H\times \Q\cdot A$.
In particular, $(H,A)$-(semi)stability is well-defined for ample $\Q$-line bundles, and $H$ or $A$ can be chosen to be very ample without changing the $(H,A)$-(semi)stability.

\begin{Lemma}\label{StabChar2d}
Let $E$ be a pure two-dimensional sheaf on a nonsingular projective surface over an algebraically closed field $k$ of characteristic zero, $H$ and $A$ two ample divisors and $F\subset E$ a nontrivial subsheaf. Then
$p_{H,A}(F)\le_0 p_{H,A}(E)$ if and only if $$\left( \mu_H(F), \frac{\chi(F)}{\rk F}, -\mu_A(F) \right)\le \left( \mu_H(E), \frac{\chi(E)}{\rk E}, -\mu_A(E) \right) \,,$$
where we consider the lexicographic ordering on $\Q^3$.
\end{Lemma}
\Proof
For a line bundle $L$ Riemann-Roch yields
$$\chi(E\otimes L)=\frac{\rk E}{2}\c_1(L)^2 + \left( \c_1(E)-\frac {\rk E} 2 K_X \right) .\c_1(L)+\chi(E) \,,$$
where $K_X$ is the canonical divisor on the surface,
and therefore
\begin{eqnarray*}
P_{H,A}(E)(m,n)
&=&\frac{\rk E}{2}(mH+nA)^2 + \left( \c_1(E)-\frac {\rk E} 2 K_X \right) .(mH+nA)+\chi(E)\\
H^2\cdot p_{H,A}(E)(m,n) \nonumber 
&=&\frac{1}{2}(mH+nA)^2 + \left( \frac{\c_1(E)}{\rk E}-\frac {1} 2 K_X \right) .(mH+nA)+\frac{\chi(E)}{\rk E}\,. 
\end{eqnarray*}
\qed

\begin{Obs}\label{StabImp}
$H$-stable $\Rightarrow$ $(H,A)$-stable $\Rightarrow$ $(H,A)$-semistable $\Rightarrow$ $H$-semistable.
\end{Obs}
This trivial observation must not be neglected: it is the reason why we can get morphisms between the corresponding moduli spaces.
Conversely, there might be $(H,A)$-stable sheaves that are not $H$-stable,
and there might be $H$-semistable sheaves that are not $(H,A)$-semistable.

$(H,A)$-(semi)stability is a generalisation of $H$-(semi)stability in the following sense:
$(H,H)$-(semi)stability is equivalent to $H$-(semi)stability,
and one has $p_H(F)=p_H(E)$ for two coherent sheaves $E$ and $F$ if and only if $p_{H,H}(F)=p_{H,H}(E)$.
In particular, everything we can prove for $(H,A)$-(semi)stability also holds for $H$-(semi)stability.

Conversely, one can generalise known facts on $H$-(semi)stability to $(H,A)$-(semi)stability.
There is always a Jordan-Hölder filtration $0=E_0\subset E_1\subset ... \subset E_\ell=E$ for an $(H,A)$-semistable sheaf $E$.
The graded object $$gr(E):=\bigoplus_{i=1}^\ell E_i/E_{i-1}$$ does not depend on the choice of the Jordan-Hölder filtration.
Two $(H,A)$-semistable sheaves $E_1$ and $E_2$ with $p_{H,A}(E_1)=p_{H,A}(E_2)$ are called Seshadri equivalent or S-equivalent if $gr(E_1)\cong gr(E_2)$.

Let $E$ be an $(H,A)$-semistable sheaf of dimension $d$ and
$0=E_0\subset E_1\subset ... \subset E_\ell=E$
a Jordan-Hölder filtration of $E$.
By observation \ref{StabImp} $E$ is in particular $H$-semistable but the factors $gr_i(E)$ are not necessarily $H$-stable.
Thus one gets a Jordan-Hölder filtration of $E$ with respect to $H$-stability by refining the given filtration.

Passing from the set of $H$-semistable sheaves to the set of $(H,A)$-semista\-ble sheaves one looses sheaves, and the S-equivalence classes become smaller. This is the reason why a moduli space for $(H,A)$-semistable sheaves parametrising $(H,A)$-polystable sheaves can partially resolve a component of a moduli space for $H$-semistable sheaves parametrising $H$-polystable sheaves.

\begin{Lemma}\label{PHAlc}
Let $f:X\to S$ be a projective morphism of noetherian schemes, $H$ and $A$ two $f$-ample invertible sheaves on $X$ and $F$ a flat family of sheaves on the fibres of $f$.
Then the polynomial $P_{H_s,A_s}(F_s)$ is locally constant as a function of $s\in S$.
\end{Lemma}
\Proof
The family $F(\ell H)$ is $S$-flat as well for all $\ell\in\N_0$, so by \cite{HL} proposition 2.1.2 the Hilbert polynomial $P_{A_s}(F_s(\ell H_s))=P_{H_s,A_s}(F_s)(\ell,\bullet)\in\Q[n]$ is locally constant as a function of $s\in S$ for all $\ell\in\N_0$.
The polynomial $P_{H_s,A_s}(F_s)$ can be regained from $P_{H_s,A_s}(F_s)(\ell,\bullet)$ for finitely many choices of $\ell$, hence it is locally constant as a function of $s\in S$ as well.
\qed

\begin{Prop}\label{StabOpen}
The following properties of coherent sheaves are open in flat families: being $(H,A)$-semistable, or $(H,A)$-stable.
\end{Prop}
\Proof
Let $f:X\to S$ be a projective morphism of noetherian schemes, $H$ and $A$ two $f$-very ample invertible sheaves on $X$ and $F$ a flat family of $d$-dimensional sheaves on the fibres of $f$ with Hilbert polynomial $P$ with respect to $H_s$ for all $s\in S$.
As we want to show the openness of certain subsets we can assume $S$ to be connected.
Furthermore, we can replace $S$ by the open subset of all $s\in S$ such that $F_s$ is $H_s$-semistable as this condition is open by \cite{HL} proposition 2.3.1, having in mind observation \ref{StabImp}.
Let $\alpha\in\N$ be the multiplicity associated to $P$.

For each $\alpha'\in\N$ with $\alpha'\le\alpha$ we consider the relative Quot scheme
$$\pi: Q(\alpha'):=\Quot_{X/S}(F,\frac{\alpha'}{\alpha} P)\to S\,,$$
see \cite{HL} section 2.2.
Let $C(\alpha')$ be the set of connected components of $Q(\alpha')$ and 
$U\in\Coh(Q(\alpha')\times_S X)$ the universal quotient family.
By lemma \ref{PHAlc} $P_{H,A}:=P_{H_s,A_s}(F_s)$ is independent of $s\in S$ and
$$P_{H,A}(C):=P_{H_{\pi(q)},A_{\pi(q)}}(U_q)$$ is independent of $q\in C$ for $C\in C(\alpha')$.
Let $p_{H,A}$ and $p_{H,A}(C)$ be the reduced polynomials associated to $P_{H,A}$ and $P_{H,A}(C)$, respectively.

Now $F_s$ is $(H_s,A_s)$-(semi)stable if and only if it is not contained in the closed union
$$\bigcup_{\alpha'=1}^\alpha \pi\left( \bigcup_{C\in C(\alpha')\;:\; p_{H,A} \lle p_{H,A}(C)} C \right)\,.$$
\qed

\section{The construction of the moduli space}\label{MPConst}

We straightforward generalise the construction in \cite{HL} chapter 4 omitting details that can be found therein.
Let $X$ be a connected projective scheme over an algebraically closed field $k$ of characteristic zero, $H$ and $A$ two ample line bundles on $X$ and $P\in\Q[x]$.
We define a functor
$\M': (Sch/k)^o\to (Sets)$
from the category opposed to the category of $k$-schemes to the category of sets as follows.
For a $k$-scheme $S$ let $\M'(S)$ be the set of all isomorphism classes of $S$-flat families of $(H,A)$-semistable sheaves on $X$ with Hilbert polynomial $P$ with respect to $H$, and for a $k$-morphism $f:S'\to S$ let
$$\M'(f): \M'(S)\to\M'(S'), [F]\mapsto [ (f\times\id_X)^* F]\,.$$
If we consider the equivalence relation $F\sim F'$ for two $F,F'\in\M'(S)$ if and only if $F\cong F'\otimes p^*L$ for some $L\in\Pic(S)$, where $p: S\times X\to S$ is the projection onto the first factor, then we get our moduli space functor as quotient functor:
$\M:=\M'/\sim$ is the moduli functor for $(H,A)$-semistable sheaves on $X$ with Hilbert polynomial $P$ with respect to $H$.

Considering only families of $(H,A)$-stable sheaves yields open subfunctors $\M'^s\subset\M'$ and $\M^s\subset\M$ as the stability condition is open in flat families, see proposition \ref{StabOpen}.
Note that as our ground field is algebraically closed, being $(H,A)$-stable is equivalent to being geometrically $(H,A)$-stable analogously to the results in \cite{HL} section 1.5.

\begin{Def}
A scheme $M$ is called a moduli space for $(H,A)$-semistable sheaves if $M$ corepresents $\M$.
We will denote $M$ by $M_{H,A}(P)$, and analogously the functors.
\end{Def}
Suppose $M$ corepresents $\M$. Then analogously to \cite{HL} lemma 4.1.2 Seshadri equivalent sheaves correspond to identical closed points in $M$.
In particular, if there is a properly $(H,A)$-semistable sheaf $F$, then $\M$ cannot be represented.

According to \cite{HL} theorem 3.3.7 the family of $H$-semistable sheaves on $X$ with Hilbert polynomial with respect to $H$ equal to $P$ is bounded. In particular, there is an integer $m$ such that any such sheaf $F$ is $m$-regular.
Let $V:=k^{\oplus P(m)}$ and $\mH:=V\otimes_k\O_X(-mH)$. Then there is a surjection
$\rho: \mH\to F\,,$
which gives a closed point
$$[\rho: \mH\to F]\in R\subset\Quot(\mH,P)\,,$$
where $\Quot(\mH,P)$ is Grothendieck's Quot scheme of quotients of $\mH$ with Hilbert polynomial $P$ on $X$, see e.g.\ \cite{HL} section 2.2,
and $R$ is the open subset of $\Quot(\mH,P)$ of all quotients $[\mH\to E]$, where $E$ is $H$-semistable and the induced map $$V\to H^0(\mH(mH))\to H^0(E(mH))$$ is an isomorphism.

Let $R^{ss}\subset R$ denote the open subscheme of those points which parametrise $(H,A)$-semistable sheaves, and
$R^s\subset R$ the open subscheme of those parametrising $(H,A)$-stable sheaves.
There is a $\Gl(V)$-action on $\Quot(\mH,P)$, and $R$, $R^{ss}$ and $R^s$ are $\Gl(V)$-invariant.

\begin{Prop}\label{Llvample}
If $\ell$ is sufficiently large then
the line bundle $$L_\ell:=\det(p_*(\mF\otimes q^*\O_X(mH+\ell A)))$$
on $\Quot(\mH,P)$ is very ample and carries a natural $\Gl(V)$-linearisation,
where $p$ and $q$ are the two projections from $\Quot(\mH,P)\times X$ to the first and second factor, respectively, and $\mF$ is the universal quotient sheaf on $\Quot(\mH,P)\times X$, see \cite{HL} section 2.2.
\end{Prop}
\Proof
If $E$ is a sheaf with Hilbert polynomial $P$ with respect to $H$ then $E(mH)$ has Hilbert polynomial 
$$P_H(E(mH))(x)= P(x+m)=:P'(x)$$
with respect to $H$.
As tensoring with a line bundle is exact, one has an isomorphism
$$\varphi: \Quot(\mH,P) \to \Quot(\mH(mH),P')$$
and $$L_\ell':=(\varphi^{-1})^* L_\ell = \det(p_*'(\mF'\otimes q'^*\O_X(\ell A)))\,,$$
where $p'$ and $q'$ are the two projections from $\Quot(\mH(mH),P')\times X$ to the first and second factor, respectively, and $\mF'$ is the universal quotient sheaf on $\Quot(\mH(mH),P')\times X$.
So we can assume without loss of generality that $m=0$.

Let $S\subset \Quot(\mH,P)$ be a connected component. The universal family $\mF$ is $\Quot(\mH,P)$-flat, hence the Hilbert polynomials $P_A(\mF_s)$ are constant on $S$ by \cite{HL} proposition 2.1.2, say $P_A(\mF_s)=P'$ for all $s\in S$.
Hence one has a closed embedding
$\psi: S \to \Quot_A(\mH,P')\,,$
where the index $A$ denotes that the Hilbert polynomial is with respect to $A$, and not to $H$ as before.
For sufficiently large $\ell$
$$L_\ell':= \det(p_*'(\mF'\otimes q'^*\O_X(\ell A)))$$
is very ample by proposition \cite{HL} 2.2.5, where $p'$ and $q'$ are the two projections from $\Quot_A(\mH,P')\times X$ to the first and second factor, respectively, and $\mF'$ is the universal quotient sheaf on $\Quot_A(\mH,P')\times X$,
and $L_\ell'$ carries a natural $\Gl(V)$-linearisation as explained in \cite{HL} section 4.3.
Thus $L_\ell|_S=\psi^* L_\ell'$ is very ample as well and carries a natural $\Gl(V)$-linearisation,
hence also $L_\ell$ itself.
\qed\\

As the center of $\Gl(V)$ is contained in the stabiliser of each point in $\Quot(\mH,P)$ we can restrict the action to $\Sl(V)$. Thus one has the notion of (semi)stable points of $\Quot(\mH,P)$ with respect to $L_\ell$ and the $\Sl(V)$-action.\\

\begin{Th}\label{StabCorres}
Suppose that $m$, and for fixed $m$ also $\ell$ are sufficiently large integers.
Then $R^{ss}=\overline R^{ss}(L_\ell)$ and $R^s=\overline R^{s}(L_\ell)$.
Moreover, the closures of the orbits of two points $[\rho_i: \mH\to F_i]$, $i=1,2$, in $R^{ss}$ intersect if and only if $gr^{JH}(F_1)\cong gr^{JH}(F_2)$. The orbit of a point $[\rho: \mH\to F]$ is closed in $R^{ss}$ if and only if $F$ is polystable.
\end{Th}
We prove this theorem later. Together with \cite{HL} theorem 4.2.10 and the analogous statement of \cite{HL} lemma 4.3.1 it yields:

\begin{Th}\label{MHAexists} 
There is a projective scheme $M_{H,A}(P)$ that universally corepresents the functor $\M_{H,A}(P)$.
Closed points in $M_{H,A}(P)$ are in bijection with Seshadri equivalence classes of $(H,A)$-semistable sheaves with Hilbert polynomial $P$. Moreover, there is an open subset $M^s_{H,A}(P)$ that universally corepresents the functor $\M^s_{H,A}(P)$.
\end{Th}
Let $[\rho: V\otimes\O_X(-mH)\to F]$ be a closed point in $\overline R$, $\lambda:\G_m\to \Sl(V)$ a one-parameter subgroup and $V=\bigoplus_{n\in\Z} V_n$ the weight space decomposition.
Define ascending filtrations on $V$ and $F$ by
$$V_{\le n}:=\bigoplus_{\nu\le n}V_\nu \mand F_{\le n}:=\rho(V_{\le n}\otimes\O_X(-mH)).$$
Then $\rho$ induces surjections $\rho_n: V_n\otimes\O_X(-mH)\to F_n:=F_{\le n}/F_{\le n-1}$. Summing up over all weights we get a closed point
$$\left[\overline\rho:=\bigoplus_{n\in\Z}\rho_n: V\otimes\O_X(-mH)\to \overline F:=\bigoplus_{n\in\Z}F_n\right]\in \Quot(\mH,P)\,.$$

\begin{Lemma}
The weight of the action of $\G_m$ via $\lambda$ on the fibre of $L_\ell$ at the point $[\overline\rho]$ is given by
$$\sum_{n\in\Z}n \chi(F_n(mH+\ell A)) = $$
$$-\frac{1}{\dim(V)}\sum_{n\in\Z}\left(\dim(V)\chi(F_{\le n}(mH+\ell A))-\dim(V_{\le n})\chi(F(mH+\ell A))  \right)\,.$$
\end{Lemma}
\Proof
This is \cite{HL} lemma 4.4.4 with minor changes due to the more general situation.
\qed

\begin{Lemma}
A closed point $[\rho: \mH\to F]\in \overline R$ is (semi)stable if and only if for all nontrivial proper linear subspaces $V'\subset V$ and the induced sheaf $F':=\rho(V'\otimes\O_X(-mH))\subset F$ the following inequality holds:
$$\dim V\cdot \chi(F'(mH+\ell A))\geg \dim(V')\cdot \chi(F(mH+\ell A))\,.$$
\end{Lemma}
\Proof
This is the generalisation of \cite{HL} lemma 4.4.5. The proof carries over literally using the replacement
$P(\bullet,\ell)\mapsto \chi(\bullet(mH+\ell A))\,.$
\qed\\

\noindent
In the following we denote $H^0(\rho(mH))^{-1}( H^0(F'(mH)))$ by $V\cap H^0(F'(mH))$.

\begin{Lemma}\label{HstPtChar}
If $\ell$ is sufficiently large, a closed point $[\rho: \mH\to F]\in \overline R$ is (semi)stable if and only if for all 
coherent subsheaves $F'\subset F$ and $V'=V\cap H^0(F'(mH))$ the following inequality holds:
$$\dim V\cdot \chi(F'(mH+z A))\geg \dim(V')\cdot \chi(F(mH+z A))$$
as polynomials in $z$.
\end{Lemma}
\Proof
This is the generalisation of \cite{HL} lemma 4.4.6. The proof carries over almost literally again.
\qed\\

\noindent Recall our choice of the ordering $\le$ on $\Q[m,n]$ explained in section \ref{Semistable}.
\begin{Lemma}\label{PHAeq}
Let $M\subset\Q[m,n]$ be a finite set of polynomials.
Then there is an $m_0\in\N$ such that for all $m'\ge m_0$ and for all $P,Q\in M$ the following conditions are equivalent:
\begin{enumerate}
\item $P\le Q$, 
\item $P(m',\bullet)\le Q(m',\bullet)$ (as polynomials in $n$) and
\item $P(m',n')\le Q(m',n')$ for some $n'\gg 0$.
\end{enumerate}
\end{Lemma}
\Proof
Straightforward.
\qed

\begin{Prop}\label{PHAfin}
Let $P\in\Q[m]$ be a polynomial.
The set $M:=\{ P_{H,A}(G) \;|\;$ $G$ is a quotient of an $H$-semistable sheaf $F$ with $P_H(F)=P$ and $p_H(F)=p_H(G) \}\subset\Q[m,n]$ is finite.
\end{Prop}
\Proof
An immediate consequence of the definition is that any $H$-destabilising quotient of an $H$-semistable sheaf is $H$-semistable as well.
Let $\alpha$ be the multiplicity associated to $P$, $\alpha'\in\N$ with $\alpha'\le\alpha$ and $\mF = (F_i)_{i\in I}$
the family of $H$-semistable sheaves with Hilbert polynomial $\alpha'p$ with respect to $H$.
By \cite{HL} theorem 3.3.7 this family is bounded. Hence the family $\mF(m' H) := (F_i(m' H))_{i\in I}$ is bounded as well for any choice of $m'\in\N$. Therefore by \cite{HL} lemma 1.7.6 the set of Hilbert polynomials $\{ P_A(F_i(m'H)) \;|\; i\in I \}$ is finite for any choice of $m'\in\N_0$.
As the polynomials $P_{H,A}(F_i)$ can be regained from $P_A(F_i(m'H))$ for finitely many choices of $m'$ the set $M_{\alpha'}:=\{ P_{H,A}(F_i) \;|\; i\in I \}$ is finite.
Altogether one has that $M$ is finite because $M\subset\bigcup_{\alpha'=1}^\alpha M_{\alpha'}$.
\qed

\pagebreak[2]
\begin{Cor}\label{HAStabNumChar}
Let $P$ be a polynomial in one variable.
Suppose that $m$, and for fixed $m$ also $\ell$ are sufficiently large integers.
Then the following conditions for an $H$-semistable sheaf $F$ with $P_H(F)=P$ are equivalent:
\begin{enumerate}
\item $F$ is $(H,A)$-(semi)stable,
\item for all nontrivial proper subsheaves $F'\subset F$ with $p_H(F')=p_H(F)$ one has
$p_{H,A}(F')$ $\geg p_{H,A}(F)$,
\item for all nontrivial proper subsheaves $F'\subset F$ with $p_H(F')=p_H(F)$ one has
$$\frac{\chi(F'(mH+z A))}{r'}\geg\frac{\chi(F(mH+z A))}{r}$$
as polynomials in $z$, where $r'$ and $r$ denotes the multiplicity of the sheaves $F'$ and $F$, and
\item for all nontrivial proper subsheaves $F'\subset F$ with $p_H(F')=p_H(F)$ one has
$$\frac{\chi(F'(mH+\ell A))}{r'}\geg\frac{\chi(F(mH+\ell A))}{r}\,.$$
\end{enumerate}
\end{Cor}
\Proof
The equivalence of conditions 1 and 2 follows immediately from the definition.
The equivalence of conditions 2-4 is established using lemma \ref{PHAeq} and proposition \ref{PHAfin}.
\qed\\

\noindent\ProofOf \textit{theorem \ref{StabCorres}}.
Let $m$ be large enough in the sense of \cite{HL} theorem 4.4.1 and of corollary \ref{HAStabNumChar} and such that any $H$-semistable sheaf with multiplicity $\rho\le r$ and Hilbert polynomial $\rho\cdot p$ with respect to $H$ is $m$-regular. 
Moreover, let $\ell$ be large enough in the sense of lemma \ref{HstPtChar}, proposition \ref{Llvample} and corollary \ref{HAStabNumChar}.

First assume that $[\rho: \mH\to F]$ is a closed point in $R$. By definition of $R$, the map $V\to H^0(F(mH))$ is an isomorphism. Let $F'\subset F$ be a subsheaf of multiplicity $0<r'<r$ and let $V'=V\cap H^0(F'(mH))$.
According to \cite{HL} theorem 4.4.1 one has either
\begin{enumerate}
\item $p_H(F')=p_H(F)$, or
\item $h^0(F'(mH))<r'\cdot p(m)$.
\end{enumerate}
In the first case $F'$ is $m$-regular, and we get $\dim(V')=h^0(F'(mH))=r'\cdot p(m)$ and therefore
\begin{eqnarray*}
&&\frac 1{rr'}\left(\dim V\cdot \chi(F'(mH+\ell A))-\dim(V')\cdot \chi(F(mH+\ell A))\right)\\
&=& \frac{\dim V}r \cdot \frac{\chi(F'(mH+\ell A))}{r'}-\frac{\dim(V')}{r'}\cdot \frac{\chi(F(mH+\ell A))}r\\
&=& \frac{\dim V}r \left(\frac{\chi(F'(mH+\ell A))}{r'}-\frac{\chi(F(mH+\ell A))}r\right).
\end{eqnarray*}
In the second case
$\dim(V)\cdot r'=rr'p(m)>h^0(F'(mH))\cdot r=\dim(V')\cdot r\,.$
These are the leading coefficients of the polynomials of lemma \ref{HstPtChar} up to a factor, so that
$$\dim V\cdot \chi(F'(mH+z A))>\dim(V')\cdot \chi(F(mH+z A))$$
as polynomials in $z$. By lemma \ref{HstPtChar} and corollary \ref{HAStabNumChar} this proves the implications
\begin{enumerate}
\item $[\rho]\in R^s \Rightarrow [\rho]\in \overline R^s(L_\ell)$,
\item $[\rho]\in R^{ss}\setminus R^s \Rightarrow [\rho]\in \overline R^{ss}(L_\ell)\setminus \overline R^s(L_\ell)$ and
\item $[\rho]\in R\setminus R^{ss} \Rightarrow [\rho]\notin \overline R^{ss}(L_\ell)$.
\end{enumerate}
Conversely, suppose that $[\rho: V\otimes\O_X(-mH)\to F]\in \overline R^{ss}(L_\ell)$.
It remains to show that $[\rho]\in R$.
By lemma \ref{HstPtChar} one has an inequality
$$\dim V\cdot \chi(F'(mH+z A))\ge \dim(V')\cdot \chi(F(mH+z A))$$
as polynomials in $z$ for any $F'\subset F$ and $V'=V\cap H^0(F'(mH))$.
Passing to the leading coefficient of the polynomials we get
$$p(m)\cdot r\cdot r'=\dim(V)\cdot r'\ge \dim(V')\cdot r\,.$$
This is the inequality (4.4) in the proof of \cite{HL} theorem 4.3.3 in chapter 4.4,
and the remaining part of this proof carries over literally.
\qed

\bibliographystyle{amsalpha}
\bibliography{my}

\end{document}